\documentclass[11pt]{amsart}

\usepackage{amssymb,latexsym}
\usepackage[mathcal]{eucal}
\usepackage{amscd}
\usepackage{amsmath}
\usepackage{graphicx}
\usepackage{amsfonts}
\usepackage{mathtools}
\usepackage{amsthm}

\usepackage{xy}
\input xy
\xyoption{all}

\newtheorem{theo}			{Theorem}[section]
\newtheorem{prop}[theo]		{Proposition}
\newtheorem{lem}[theo]		{Lemma}
\newtheorem{cor}[theo]		{Corollary}
\newtheorem{df}[theo]		{Definition}
\newtheorem{ex}[theo]		{Example}
\newtheorem*{ques}			{Question}

\def \ol {\overline}
\def \N {\mathbb N}
\def \Z {\mathbb Z}
\def \R {\mathbb R}
\def \a {\alpha }
\def \b {\beta}
\def \c {\mathfrak{c}}

\long\def\symbolfootnote[#1]#2{\begingroup%
\def\thefootnote{\fnsymbol{footnote}}\footnote[#1]{#2}\endgroup} 
\newcommand		{\D}		{\mathcal{D}}
\newcommand           {\0}		{\varnothing}
\renewcommand	{\o}         	{\circ}
\newcommand		{\nd}		{\noindent}
\newcommand		{\Trans}     {\textit{Trans}_f}
\newcommand		{\Transz}	{\textit{Trans}_s^0}
\newcommand		{\Iso}		{\textit{Iso}_X}
\newcommand		{\nimplies}{\ \not \!\!\!\!\implies}
\newcommand		{\x}		{\times}
\newcommand		{\rest}	{\upharpoonright}
\newcommand		{\iso}		{\cong}
\newcommand		{\ang}[1]    {\langle #1 \rangle}
\renewcommand	{\-}         	{^{-1}}
\newcommand		{\less}        {\setminus}
\newcommand		{\comment}[1]              {}
\newcommand		{\mr}         	{\mathrm}
\newcommand		{\Union}    {\bigcup}
\renewcommand	{\:}		{\vcentcolon}
\mathchardef\mhyphen="2D

\numberwithin{equation}{section}

\begin{document}

\bibliographystyle{alpha}

\title{Conceptions of Topological Transitivity}

\author{Ethan Akin}
\address{
Mathematics Department, The City College, 137 Street and Convent Avenue, 
New York, NY 10031, USA; (212) 650-5136}
\email{ethanakin@earthlink.net} 

\author{Jeffrey D. Carlson}
\address{Mathematics Department, 
Tufts University, Medford, MA 02155}

\date{\today}

\begin{abstract}
There are several different common
definitions of a property in topological dynamics
called ``topological transitivity,''
and it is part of the folklore of dynamical systems that under reasonable hypotheses,
they are equivalent. 
Various equivalences are proved in different places, but the full story
is difficult to find. 
This note provides a complete description of the relationships among the different
properties.
\end{abstract}

\subjclass{Primary 37B05; Secondary 37B99, 54H20}

\keywords{Topological transitivity; transitive points; minimal sets; isolated points}

\maketitle

\vspace{1cm}

\section{Introduction}
A discrete-time \emph{dynamical system} $(X,f)$ is a continuous map $f$ on a 
nonempty topological space $X$, i.e.
$f\colon X \to X$.
The dynamics is obtained by iterating the map. 

Let $\Z$ and $\N$ be the additive group of integers and the semigroup of
nonnegative integers, respectively. 
A dynamical system $(X,f)$ induces an action of $\N$ on $X$ by $n \mapsto f^n$,
where $f^0(x) = x$ and $f^{n+1}(x) = f\big(f^n(x)\big)$ for all $x \in X$ and $n \in \N$. 
A dynamical system $(X,f)$ with $f$ a homeomorphism 
is called \emph{invertible} or reversible with \emph{inverse} $(X,f^{-1})$. 
When $f$ is a homeomorphism we obtain an action of $\Z$ on $X$.


For $A \subset X$ and $k \in \Z$ we denote by $f^{k}(A)$
the image of $A$ under $f^k$ when $k \geq 0$
and the preimage under $f^{|k|}$ when $k < 0$. 
In the case of a single point, $x$, 
we will write $f^{-k}(x)$ for the set $f^{-k}(\{ x \})$ for $k > 0$ 
and let context determine in the invertible case 
whether $f^{-k}(x)$ refers to the singleton set or the point contained therein.

We say  $A \subset X$  is \emph{$+$invariant}
when $f(A) \subset A$, or equivalently $A \subset f^{-1}(A)$,
and $A$ is \emph{$-$invariant} 
when $f^{-1}(A) \subset A$.
Clearly, $A$ is $+$invariant iff its complement is $-$invariant. 
We call $A$ \emph{invariant} when $f(A) = A$.
When $A$ is $+$invariant the restricted dynamical system $(A, f \rest A)$  
is called a \emph{subsystem} of $(X,f)$.

An action of a group $G$ on a set $S$ is called \emph{transitive} 
when there is a single orbit;
that is, for any $x, y \in S$ there exists $g \in G$ such that $g x = y$, 
or equivalently, $S$ 
does not contain a proper invariant subset.
For a dynamical system $(X,f)$,
notions of ``topological transitivity'' are obtained by
replacing the original points by arbitrarily close approximations.
Differing notions of approximation yield a number of distinct properties 
that have historically been used as definitions of topological transitivity.
The relationship between the various properties is part of the folklore of the subject,
but it is hard to find a reference where the issues are sorted out.
We take this opportunity to do so. Some of the details appear to be new.

\medskip


For $x \in X$ the associated \emph{forward orbit} is 
the smallest $+$invariant set containing $x$, 
\[
O(x) \: = \{ f^k(x) : k \in \N \}.
\]
The smallest $+$invariant set containing $A$ is 
$O(A) \: = \Union_{x \in A} O(x)$ for $A \subset X$;
$A$ is $+$invariant iff $A = O(A)$.
For any subset $A \subset X$, 
the smallest $-$invariant set containing $A$ is
\[
O_-(A) \: = \Union_{k \in \N} f^{-k}(A),
\]
and $A$ is $-$invariant iff $A = O_-(A)$.
We set $O_-(x) = O_-\big(\{x\}\big)$.
We will also write 
\[O_{\pm}(A) \: = O(A) \cup O_-(A) = \Union_{k \in \Z} f^k(A).\]
This set is $+$invariant, but not necessarily $-$invariant unless $f$ is injective.
It is invariant if $f$ is bijective.
The case $O_{\pm}(x) = O_\pm\big(\{x\}\big)$ will be important later. 
The \emph{omega-limit set} of $x$,
denoted $\omega f(x)$, is the set of
limit points of the forward orbit $O(x)$:
\[
\omega f(x)  \: = 
\bigcap_{n \geq 0} \overline{\{ f^k(x) : k \geq n \}} = 
\bigcap_{n \geq 0} \overline{O\big(f^n(x)\big)}.\]

A bi-infinite sequence $\ang{ x_k : k \in \Z }$ 
is an \emph{orbit sequence} when $f(x_k) = x_{k+1}$ for all $k$;
the set $\{x_k : k \in \Z \}$ of its elements is called an \emph{orbit}.
We will also call a sequence $\ang{ x_k : k \geq n }$ an orbit sequence if
$f(x_k) = x_{k+1}$ for $k \geq n$ and $f^{-1}(x_n) = \0$;
the set of elements of this sequence is $O(x_n)$;
its set of elements is again called an orbit. 
Thus an orbit sequence is $+$invariant and cannot be extended backward.
$O_\pm(x)$ is the union of all orbits with $x_0 = x$.
Because of our two-pronged definition there is for every
$x \in X$ at least one orbit sequence with $x_0 = x$. 
When $f$ is bijective, this orbit sequence is unique and we call $O_{\pm}(x)$
\emph{the orbit of} $x$.


For $A, B \subset X$ note that
$f^k(A) \cap B \neq \0 \iff A \cap f^{-k}(B) \neq \0$ for $k \in \N$ because
each says there exist $x \in A$ and $y \in B$ such that $ y = f^k(x)$. 
Define the \emph{hitting time sets}
\begin{align*}\label{1}
N(A,B)  & \: =  \{ k \in \Z : A \cap f^{-k}(B) \neq \0 \}, \\
N_+(A,B)  & \: =  N(A,B) \cap \N.
\end{align*}
Thus for $k \geq 0$ we have $k \in N(A,B) \iff k \in N_+(A,B)$, 
and for $k \leq 0$ we have $k \in N(A,B) \iff -k \in N_+(B,A)$.
It follows that
\begin{equation}\label{neg}
N(A,B) = N_+(A,B) \cup -N_+(B,A).
\end{equation}
It also follows from our definitions above that 
\begin{equation}\label{NUVOV}
N_+(A,B) \neq \0 
\iff O(A) \cap B \neq \0
\iff A \cap O_-(B) \neq \0;
\end{equation}
\begin{equation}\label{NUVOpmV}
N(A,B) \neq \0 
\iff O_\pm(A) \cap B \neq \0
\iff A \cap O_\pm(B) \neq \0.
\end{equation}

As the condition will recur constantly, a set $U$ will be called 
\emph{opene} when it is open and nonempty. 
A sequence will be said to be \emph{dense} whenever its associated set is.

A collection $\D$ of opene subsets of a space $X$ will be 
called a \emph{density basis} 
when, given any $A \subset X$,
if for all $U \in \D$ we have $A \cap U \neq \0$,
then $A$ is dense in $X$. 
Equivalently, $\D$ is a density basis if $X$ is the only closed set which meets every
element of $\D$ 
(see Proposition \ref{denprop1} in the appendix on density bases).
The property of admitting a countable density basis will be important 
in finding dense orbits.

\begin{df}Let $(X,f)$ be a dynamical system.
To describe topological transitivity, we label seven possible properties of $(X,f)$.

\begin{itemize}
\item[(IN)]
$X$ is not the union of two proper, closed, $+$invariant subsets.

\item[(TT)]
For every pair $U,V$ of opene subsets of $X$, 
the set $N(U,V)$ is nonempty.
In this case, we say $(X,f)$ is \emph{topologically transitive}.

\item[(TT$_+$)]
For every pair $U,V$ of opene subsets of $X$, 
the set $N_+(U,V)$ is nonempty.

\item[(TT$_{++}$)]
For every pair $U,V$ of opene subsets of $X$, 
the set $N_+(U,V)$ is infinite.

\item[(DO)]
There exists an orbit sequence $\ang{ x_k : k \in \Z }$ or $\ang{ x_k : k \geq n }$
dense in $X$.

\item[(DO$_+$)] There exists a point $x \in X$ 
with forward orbit $O(x)$ dense in $X$. 
In this case, we say $(X,f)$ is \emph{point transitive}, 
and call $x$ a \emph{transitive point}.

\item[(DO$_{++}$)] 
There exists $x \in X$ such that the omega-limit set $\omega f(x) = X$.
\end{itemize}
The set of transitive points of $X$ is labelled $\Trans$.
Thus when $\Trans$ is nonempty, the system $(X,f)$ is point transitive.
When $\Trans = X$, the system $(X,f)$ is called \emph{minimal}.

\end{df}

It is clear that
\begin{equation}\label{impl1}
\vcenter{\xymatrix@R=0pt{
\mr{DO}_{++}  	\ar@{=>}[r] &
\mr{DO}_+  	\ar@{=>}[r] &
\mr{DO};					 \\
\mr{TT}_{++}	\ar@{=>}[r] &
\mr{TT}_+ 		\ar@{=>}[r] &
\mr{TT}.
}
}
\end{equation}

We can add two vertical implications:
\begin{prop}\label{DOTT}
Let $(X,f)$ be a dynamical system.
\begin{enumerate}
\item[(a)] If $X$ has a dense orbit sequence, $(X,f)$ is topologically transitive, i.e.
$\mr{DO} \implies \mr{TT}$.
\item[(b)] If there exists $x \in X$ so that $\omega f(x) = X$,
then for every pair $U,V$ of nonempty, open subsets of $X$, $N_+(U,V)$ is infinite,
i.e. $\mr{DO}_{++} \implies \mr{TT}_{++}$.
\end{enumerate}
\end{prop}

\begin{proof}
(a) Suppose $(X,f)$ satisfies $\mr{DO}$.
    Let opene $U,V \subset  X$
    and a dense orbit sequence $\ang{ x_k }$ of $f$ be given.
    By assumption there are $n,m \in \mathbb{Z}$
    such that $x_n \in U$ and $x_m \in V$.
    Since $x_n \in f^{n-m}(\{x_m\})$,
    we have $x_n \in U \cap f^{n-m}(V)$, and so $(X,f)$ satisfies $\mr{TT}$.
    
(b) Suppose $(X,f)$ satisfies $\mr{DO}_{++}$.
Let $x \in X$ be such that $\omega f(x) = X$ and let opene $U,V \subset X$
be given.
Since $U,V \subset X = \ol{O\big(f^n(x)\big)}$ for each $n \in \N$,
there are $j \in \N$ such that $f^j(x) \in U$ and infinitely many $k \geq j$
such that $f^k(x) \in V$.
Thus $f^j(x) \in U \cap f^{-(k-j)}(V) $ for infinitely
many $k \geq j$. Hence $(X,f)$ satisfies $\mr{TT}_{++}$.
\end{proof}

The $\mr{TT}$ condition is equivalent to $\textrm{IN}$:
\begin{prop}\label{INTT}
For a dynamical system $(X,f)$ the following conditions are equivalent:
\begin{itemize}
\item[(i)] $(X,f)$ is topologically transitive. $\mr{(TT)}$

\item[(ii)] $X$ does not contain two disjoint, opene, $-$invariant subsets.

\item[(iii)] $X$ is not the union of two proper, closed, $+$invariant subsets. 
$\mr{(IN)}$
\end{itemize}
\end{prop}

\begin{proof} (i) $\implies$ (ii):
If $U_1$ and $U_2$ are disjoint opene, $-$invariant subsets of $X$
then $N_+(U_1,U_2) \cup -N_+(U_2,U_1) = N(U_1,U_2)$ is empty 
and so $(X,f)$ is not topologically transitive.

(ii) $\implies$ (i): If $U_1$ and $U_2$ are opene subsets of $X$,
then the opene sets $O_-(U_1)$ and $O_-(U_2)$ meet by (ii),
and so there exist $x \in X$ and $k_i \geq 0$
such that $f^{k_i}(x) \in U_j$ for $j = 1,2$.
Without loss of generality we can assume $k_1 \leq k_2$.
Then $y \: = f^{k_1}(x) \in U_1$ and $f^{k_2 - k_1}(y) \in U_2$,
so $k_2 - k_1 \in N_+(U_1,U_2) \neq \0$.

(iii) $\iff$ (ii):
$A_1, A_2$ are closed, proper, $+$invariant subsets with union $X$ iff
the complements $X \setminus A_1, X \setminus A_2$ 
are opene, $-$invariant subsets which are disjoint.
\end{proof}

We thus have the following implications for any dynamical system:

\medskip

\centerline{
\xymatrix{
\mr{DO}_{++}  	\ar@{=>}[r] 	\ar@{=>}[d] 	&
\mr{DO}_+  	\ar@{=>}[r]  				&
\mr{DO} 		\ar@{=>}[d] 				&	\\
\mr{TT}_{++} 	\ar@{=>}[r] 				&
\mr{TT}_+		\ar@{=>}[r] 				&
\mr{TT} 		\ar@{<=>}[r] 				&
\mr{IN}.
}
}

\bigskip

When $X$ is a Hausdorff space, the difficulties in reversing the horizontal arrows 
are associated with the occurrence of isolated points.
A point $x $ in a Hausdorff space $ X$ is \emph{isolated} when 
the singleton containing it is open, 
and hence \emph{clopen} (= closed and open).
When there are no isolated points, or equivalently, every opene set is infinite, 
then $X$ is called \emph{perfect}.

Obtaining the $\mr{DO}$ conditions from $\mr{TT}$ 
requires some strong topological hypotheses. 
The appropriate assumptions are that 
$X$ is \emph{second countable}, meaning that $X$ has a countable base, 
and that $X$ is \emph{non-meager} (or \emph{of second Baire category}), 
meaning it is not the union of a countable family of nowhere dense subsets. 
For example, if $X$ is a locally compact, separable metric space, 
these conditions hold.

\begin{theo}\label{theo} Let $(X,f)$ be a dynamical system with $X$ a Hausdorff space.
\begin{enumerate}
\item[(a)] If $X$ is perfect,
then $\mr{TT}_{++}$ is equivalent to $\mr{TT}_+$, $\mr{TT}$, and $\mr{IN}$.

\item[(b)] If $X$ is perfect,
then $\mr{DO}_{++}$ and $\mr{DO}_+$  are equivalent.

\item[(c)] If $X$ is second countable and non-meager, then $\mr{TT}_+$ implies $\mr{DO}_+$.
\end{enumerate}
\end{theo}
Thus if $X$ is perfect, Hausdorff, 
second countable, and non-meager,
then all seven conditions are equivalent.

\medskip

In Section 4, we prove the theorem, showing that if $X$ is
perfect, second countable, and non-meager, then the weakest condition $\mr{TT}$
implies all the others.
In Section 5, 
we describe in detail what can happen when the space has isolated points.
In Section 6,
we consider minimality.
In Section 7, we provide examples
which show that various other possible implications fail under weaker hypotheses.
The first appendix deals with the apparatus of density bases and
the second describes two implications that survive 
even if the space $X$ is not Hausdorff.

\section{The literature}\label{earlier}
We briefly
describe some
implications
that have been proven in the literature.
\begin{itemize}
\item In \cite{S}, it is shown that $\mr{DO}_+ \implies \mr{TT}_+$
for a metric space with no isolated points and
$\mr{TT}_+ \implies \mr{DO}_+$ for separable, non-meager metric spaces.
\item \cite {KS2} states, but does not prove,
that in a compact metric space,
sixteen conditions, including $\mr{TT}_+$, 
are equivalent and imply $\mr{DO}_+$,
and that for a perfect space, $\mr{DO}_+ \implies \mr{TT}_+$.
The sources they cite focus on compact metric spaces.
\item \cite{BC} states without proof that
$\mr{DO} \iff \mr{TT} \iff \mr{TT}_+$ 
for homeomorphisms of compact metric spaces.
\item \cite{W} shows that $\mr{DO}_+ \iff \mr{TT}_+$ 
for continuous, surjective maps on second countable Baire spaces 
(his hypotheses can be weakened to these).
\item \cite{V} shows $\mr{DO} \iff \mr{TT}$ for compact metric spaces.
\item In \cite{KH}, the equivalence of 
$\mr{DO}$, $\mr{DO}_+$, $\mr{TT}$, and $\mr{TT}_+$ is stated,
but the proof given is incomplete.
In their later work \cite{HK2}, the problem is partially fixed.
\end{itemize}

\section{Lemmas}



\begin{prop}\label{InvOpene} 
Let $(X,f)$ be a dynamical system.
\begin{enumerate}

\item[(a)] The system $(X,f)$ satisfies $\mr{TT}_+$
iff every $-$invariant opene subset of $X$ is dense.

\item[(b)] If the system $(X,f)$ satisfies $\mr{TT}_+$,
then every $+$invariant opene subset of $X$ is dense. 
The converse holds if $f$ is open.

\comment{
\item[(c)] Assume $f$ is injective. 
Then $(X,f)$ is topologically transitive
iff every subset that is both $+$invariant and $-$invariant and has nonempty interior
is dense.
}															

\item[(c)] If $f$ is a homeomorphism, 
then $(X,f)$ is topologically transitive
iff every invariant opene subset of $X$ is dense.
\end{enumerate}
\end{prop}

\begin{proof}

(a) 
By (\ref{NUVOV}), 
$\mr{TT}_+$ holds just if $O_-(V)$ is dense for each opene $V$. 
Since each $O_-(V)$ is opene and $-$invariant and 
every $-$invariant opene $W$ is $O_-(W)$,
the result follows.

\comment{
If every $-$invariant opene subset of $X$ is dense then $(X,f)$ is
topologically transitive by (ii) $\implies$ (i)
of Proposition \ref{INTT}.  
Now assume $(X,f)$ is topologically transitive,
$X$ is perfect and $V$ is an opene, $-$invariant open set.
For any opene subset $U$, $N_+(U,V)$ is nonempty by $\mr{TT}_+$ 
(Proposition \ref{prop1}). 
That is, $U \cap f^{-k}(V) \neq \0$ for some $k \in {\mathbb N}$. 
Since $V$ is $-$invariant, $f^{-k}(V) \subset V$ and so
$U$ meets $V$.
As $U$ was arbitrary, $V$ is dense.
}													

(b)
By (\ref{NUVOV}), 
$\mr{TT}_+$ holds just if $O(U)$ is dense for each opene $U$.
If $\mr{TT}_+$ holds and $U \subset X$ is a $+$invariant opene set,
then $U$ is of the form $O(U)$, and hence is dense. 
If all $+$invariant opene sets are dense and
$f$ is open, then $O(U)$ is dense for every opene $U$,
so $\mr{TT}_+$ holds.

\comment{
(c)

Recall that sets $O_\pm(A)$ are both $+$invariant and $-$invariant because $f$ 
is injective.
Assume $\mr{TT}$ and let $A \subset X$ be $+$invariant and $-$invariant 
with nonempty interior. Then each opene $U$ meets
$O_\pm(A^\o) \subset O_\pm(A) = A$, so $A$ is dense.
On the other hand, given an opene $V$, the set $O_\pm(V)$ 
is $+$invariant and $-$invariant and contains $V$,
so if each such set is dense, it follows each opene $U$ meets it 
and so $(X,f)$ satisfies $\mr{TT}$.
}															

(c) 
By (\ref{NUVOpmV}),
$\mr{TT}$ holds just if $O_\pm(U)$ is dense for each opene $U$.
Since $f$ is a homeomorphism these sets $O_\pm(U)$ are all invariant and opene; 
and if an opene set $U$ is invariant, we have $U = f(U) = f\-(U)$ by bijectivity,
and therefore $O_\pm(U) = U$, so these are the only opene invariant sets.
The result follows.
\end{proof}

\begin{prop}
Let $(X,f)$ be a dynamical system and $x \in X$.
Then $\omega f(x) = X$ if and only if $N_+(\{x\},U)$ is infinite for each opene $U$.
\end{prop}
\begin{proof}
$\omega f(x) = \bigcap_{k \in \N} \ol{O\big(f^k(x)\big)} = X$
if and only if each $O\big(f^k(x)\big)$ is dense in $X$.
This happens if and only if each meets each opene $U \subset X$,
which in turn happens if and only if for all $k> 0$, there exists 
$n > k$ such that $f^k(x) \in U$. 
But that is the same as saying $N_+(\{x\},U)$ is infinite for each $U$.
\end{proof}

\begin{cor}\label{prop3b}
Let $(X,f)$ be a dynamical system.
\begin{enumerate}
\item[(a)] $\Trans$ is a $-$invariant subset of $X$.

\item[(b)] If $X$ is a perfect $T_1$ space and $\Trans$ is nonempty,
then $\Trans$ is a $+$invariant, dense subset of $X$.
\end{enumerate}
\end{cor}
\begin{proof}
(a)
If $f(x) \in \Trans$, then $O\big(f(x)\big)$ is dense.
Thus $O(x) \supset O\big(f(x)\big)$ is dense,
so $x \in \Trans$.

(b)
\comment{
Let $x \in \Trans$.
Since $X$ is perfect,   Proposition \ref{prop3a} implies
$\omega f(x) = X$.
Clearly for any point $y$ we have $\omega f\big(f(y)\big) = \omega f(y)$,
so $\omega f\big(f(x)\big) = X$, and \emph{a fortiori}
$O\big(f(x)\big)$ is dense as well.
Thus $f(x) \in \Trans$ and
it follows that $\Trans$ is $+$invariant.
Since $\Trans$ contains the dense set $O(x)$, it is itself dense.}			
Let $x \in \Trans$. 
Since $O(x)$ is dense and $X$ is perfect $T_1$,
$O(x) \less \{x\} = O\big(f(x)\big)$ is also dense\footnote{\ 
See the footnote in the second appendix.%
},
so $f(x) \in \Trans$.
Thus $\Trans$ is $+$invariant; 
and also $O(x) \subset \Trans$, 
so $\Trans$ is dense.
\end{proof}

\section{Equivalences}

\noindent\emph{N.B. 
From now until the second appendix we will assume 
that the state space $X$ of our dynamical system $(X,f)$ is a Hausdorff space.}

\medskip

\comment{
Recall that for $A, B \subset X$, we have defined
\begin{align}\label{1again}
N(A,B)  & =  \{ k \in \Z : A \cap f^{-k}(B) \neq \0 \}; \\
N_+(A,B)  & =  N(A,B) \cap \N.
\end{align}
Note that
\begin{equation}\label{2}
k \in N(A,B) \iff -k \in N(B,A).
\end{equation}
In particular, $N(A,A) \subset \Z$ is symmetric about zero.

We are using condition $\mr{TT}$ as our definition of topological transitivity.  
That is, we call $(X,f)$ \emph{topologically transitive} 
if for all opene subsets $U, V \subset X$,
$N(U,V) \neq \0$.

From (\ref{2}) it is clear that $(X,f)$ is topologically transitive iff for every
opene $U, V \subset X$
\begin{equation}\label{3}
N_+(U,V) \cup N_+(V,U) \neq \0.
\end{equation}
Notice that if $f$ is bijective then $-N(A,B)$ is the hitting time set  for $f^{-1}$. 
This implies that a homeomorphism $f$ on $X$
is topologically transitive iff its inverse $f^{-1}$ is topologically transitive.
Similarly, if $f$ is invertible and $\ang{ x_n }$ is a dense orbit sequence for $(X,f)$,
then $\ \ang{ x_{-n} }$ is a dense orbit sequence for $(X,f^{-1})$.
}

 The key assumption needed to prove
 $\mr{TT} \implies \mr{TT}_+ \implies \mr{TT}_{++}$
 is that $X$ is perfect.


\begin{lem}\label{NUU}
If $X$ is perfect and $(X,f)$ is topologically transitive,
then for any opene $U \subset X$, the set $N_+(U,U)$ is infinite.
\end{lem}
\begin{proof}
We will by induction find a decreasing sequence
of opene sets $U_n \subset U$ and an strictly increasing sequence $k_n \in \N$
such that $f^{k_n}(U_n) \subset U$.

For the initial step, $U$ contains two distinct points because $X$ is perfect. 
As $X$ is Hausdorff, these have disjoint neighborhoods $V_1, W_1 \subset U$.
By $\mr{TT}$, there is $k_1 \in N_+(V_1,W_1) \cup N_+(W_1,V_1)$,
so that $U_1 \: = U \cap f^{-k_1}(U)$ is opene;
evidently $f^{k_1}(U_1) \subset U$. 
Because $V_1$ and $W_1$ are disjoint, $k_1 > 0$.

For the induction step, suppose we have $U_n \subset U$
and $k_n  > 0$
such that $f^{k_n}(U_n) \subset U$.
Applying the initial step to $U_n$ 
gives us an opene $U_{n+1} \subset U_n$ and a $j_{n+1} > 0$
such that $f^{j_{n+1}}(U_{n+1}) \subset U_n$;
then setting $k_{n+1} = k_n + j_{n+1} $,
we have $f^{k_{n+1}}(U_{n+1}) \subset f^{k_n}(U_n) \subset U$ 
and $k_{n+1} > k_n$.

Clearly, $k_{n} \in N_+(U,U)$ for all $n$.
\end{proof}

\begin{prop}\label{prop1}
If $X$ is perfect and $(X,f)$ is topologically transitive,
then
for all opene $U,V$, the set $N_+(U,V)$ is infinite; i.e.
$\mr{TT} \implies \mr{TT}_{++}$.
\end{prop}
\begin{proof}

$\mr{TT}  \implies \mr{TT}_+$:
Let $U,V \subset X$ be opene. 
$\mr{TT}$ says $N(U,V) \neq \0$, so by (\ref{neg})
it will be enough to show that $N_+(U,V) \neq \0$ iff $N_+(V,U) \neq \0$. 
To that end, assume without loss of generality there exists an $n \in N_+(U,V)$;
we will find an element of $N_+(V,U)$.

Now $W \: = U \cap f^{-n}(V)$ is opene, 
so by Lemma \ref{NUU}, $N_+(W,W)$ is infinite, and 
thus there exists $k > n$ such that $f^k(W) \cap W \neq \0$.
Since 
\[
f^k(W) = 
f^k\big(U \cap f^{-n}(V)\big) \subset 
f^k\big(f^{-n}(V)\big) = 
f^{k-n}(V)
\]
and $W \subset U$,
it follows that $f^{k-n}(V) \cap U \neq \0$, so that $k - n \in N_+(V,U)$.

$\mr{TT}_+ \implies \mr{TT}_{++}$:
We show that for any opene $U,V \subset X$, 
if $N_+(U,V) \neq \0$, then $N_+(U,V)$ is infinite.

Suppose $n \in N_+(U,V)$, so $W \: = U \cap f^{-n}(V)$ is opene.
For any $k \in N_+(W,W)$, we have $f^k(W) \cap W \neq \0$, 
and since $f^n(W) \subset V$, 
we then have $\0 \neq f^{k+n}(W) \cap f^n(W) \subset
f^{k+n}(W) \cap V$.
Because $W \subset U$, we have $k+n \in N_+(U,V)$.
Thus $N_+(W,W) + n \subset N_+(U,V)$. 
By Lemma \ref{NUU}, $N_+(U,V)$ is infinite.
\comment{
By Lemma \ref{NUU} there are arbitrarily large $k \in N_+(W,W)$,
so that $f^{k}(W) \cap W \neq \0$.
Since $W \subset f^{-m}(V)$,
we then have $f^{k+m}(W) \cap V \neq \0$.
Thus $N_+(W,W) + k \subset N_+(U,V)$, making the latter infinite.
Since we also have infinitely many $k > m$
such that $f^k(W) \cap W \neq \0$,
the sets $U \supseteq W$ and $f^{k-m}(V) \supseteq f^k(f^{-m}(V)) \supseteq f^k(W)$
meet for infinitely many $k-m$, so
$N_+(V,U) \supseteq [N_+(W,W) - k] \cap \N$ is infinite.}
\end{proof}

\begin{lem}\label{lem1xxx} If $(X,f)$ is topologically transitive and $A$ is a  
$+$invariant subset of $X$
then the interior of $f^{-1}(A) \setminus A$
either is empty or consists of a single isolated point.
\end{lem}

\begin{proof}
Suppose that the interior $U$ of  $f^{-1}(A) \setminus A$ contains at least two points.
Then it contains
disjoint opene subsets $U_1, U_2$. Since $f^k(U_j) \subset A$ for all $k > 0$, 
$N_+(U_1,U_2) \cup N_+(U_2,U_1) = \0$.
Hence $(X,f)$ is not topologically transitive.
If the open set $U$
is a single point then the point is isolated.
\end{proof}

\begin{cor}\label{cor1xxx} 
If $(X,f)$ is topologically transitive then either $f(X)$ is dense in $X$ or 
$X \setminus \overline{f(X)}$ consists of a single isolated point. 
In particular, if $X$ is perfect then $f(X)$ is dense.
\end{cor}

\begin{proof}
Apply the lemma with $A$  the image $f(X)$ so that $f^{-1}(A) = X$.
\end{proof}

\begin{prop}\label{prop3xxxa} 
If 
$(X,f)$ is topologically transitive and 
$A$ is a proper, closed, $+$invariant subset of $X$ 
with nonempty interior $U = A^{\circ}$, 
then
\[
f^{-1}(U) \setminus A = 
(f^{-1}(A) \setminus A)^{\circ} = 
f^{-1}(A)^\o \less U =
f^{-1}(U) \setminus U.
\]
is an isolated singleton.
Furthermore, $A \setminus U$ is $+$invariant.
\end{prop}

\begin{proof}
Since $A$ is proper, closed, and $+$invariant, 
the complement $X \setminus A$ is opene and $-$invariant. 
If $f^{-1}(U) \subset A$ then since $U$ is the interior of $A$, 
$f^{-1}(U) \subset U$ and so
$U$ would be an opene, $-$invariant subset disjoint from $X \setminus A$.  
Since $(X,f)$ is transitive, this contradicts Proposition \ref{INTT}.  
It follows that  $f^{-1}(U) \setminus A$ is nonempty. 
As it is open, it is a subset of $\big(f^{-1}(A) \setminus A\big)^{\circ}$. 
By Lemma \ref{lem1xxx}, the latter set is an open singleton $\{ x \}$; 
hence these two sets are equal.
Since $\big(f^{-1}(A) \setminus A\big)^{\circ} =
f^{-1}(A)^\o \setminus A = \{x\}$, 
we see the open set $f^{-1}(A)^\o \setminus \{ x \} \subset A$, 
so $f^{-1}(A)^\o \setminus \{ x \} \subset A^{\circ} = U$.
Hence $f^{-1}(A)^\o \setminus U = \{ x \}$. 
Since $f^{-1}(U) \less A \subset f^{-1}(U) \less U \subset f^{-1}(A)^\o \less U$,
we also see $f^{-1}(U)\less U = \{x\}$.

Finally, let $y \in A \setminus U$. 
Then $f(y) \in A$. 
Since $x \notin A$, 
we have $y \neq x$,
so $y \notin f^{-1}(U) \less U$. 
Thus $f(y) \notin U$, 
so $f(y) \in A \less U$.
\end{proof}

\comment{
\begin{cor}\label{cor3xxxa} If 
$X$ is perfect and $(X,f)$ is topologically transitive, 
then any opene, $+$invariant subset is dense in $X$.
\end{cor}
\begin{proof}
If $A$ is the closure of an opene, $+$invariant subset, then it is a closed,
$+$invariant subset with nonempty interior. 
If $A$ were proper, by Proposition \ref{prop3xxxa} 
the member of the singleton $f^{-1}(A^\o) \less A$ would be an isolated point.
\end{proof}
}

\begin{prop}\label{prop3a}
Let $X$ be perfect, $(X,f)$ point transitive, and $x \in \Trans$.
Then $\omega f(x) = X$,
so $\mr{DO}_+ \implies \mr{DO}_{++}$.
\end{prop}

\begin{proof}
$f$ satisfies $\mr{TT}$ by (\ref{impl1}) and Proposition \ref{DOTT}(a): 
$\mr{DO}_+ \implies \mr{DO} \implies \mr{TT}$.
Because $X$ is perfect,
Corollary \ref{cor1xxx} implies that $f(X)$ is dense in $X$.
Thus the preimage $f^{-1}(U)$ of each opene $U$ is opene;
by induction, $f^{-k}(U)$ is opene for each $k \geq 0$. 
Since $x \in \Trans$, we have $O(x) \cap f^{-k}(U) \neq \0$,
so that $U$ meets $f^k\big(O(x)\big) = O\big(f^k(x)\big)$.
As $U$ was arbitrary, 
each $O\big(f^k(x)\big)$ is dense, so that $\omega f(x) = X$.
\end{proof}




\begin{prop} \label{prop5} 
If $(X,f)$ is a dynamical system satisfying $\mr{TT}_+$, 
with $X$ a non-meager space admitting a countable density basis, 
then $\Trans$ is a $G_\delta$ subset.
Thus $\mr{TT}_+ \implies \mr{DO}_+$ under these hypotheses.
If $X$ is perfect, $\Trans$ is a dense $G_\delta$.
\end{prop}

\begin{proof} 

Let $x \in X$ and let $\{U_n\}_{n \in \N}$ be a density basis $X$.
By the definitions, $x \in \Trans$ iff $O(x)$ is dense iff 
$O(x)$ meets each opene $U$ iff $O(x)$ meets each $U_n$.
But from (\ref{NUVOV}), 
$O(x)$ meets $U_n$ just if $x \in O_-(U_n)$, 
so that $\Trans = \bigcap_{n \in \N} O_-(U_n)$. 

By Proposition \ref{InvOpene}(a), each $O_-(U_n)$ is dense. 
Thus we obtain $\Trans$ as the 
intersection of a countable family of dense open sets.
It is nonempty because $X$ is not meager. 
If $X$ is perfect, $\Trans$ is dense by Corollary \ref{prop3b}. 
\end{proof}



Beyond compact spaces,
we note that the broadest dynamically 
useful class to which Proposition \ref{prop5} applies
is that of Polish spaces.
A space is \emph{Polish} when it admits a complete, 
separable metric compatible with the topology.
Separable metric spaces are second countable 
and by the Baire Category Theorem any countable intersection of
dense open subsets of a completely metrizable space is dense in the space.
In particular, a Polish space is non-meager.
A subset $A$ of a complete metric space $X$ admits a complete metric
iff $A$ is a $G_{\delta}$ subset, i.e., a countable intersection of open sets
(for this result as an exercise with hints and for references, see \cite{K}, p. 207).
Thus a $G_{\delta}$ subset of a Polish space is Polish.
If $X$ is a second countable, locally compact space,
then its one point compactification is a second countable, compact space
and so $X$ is a $G_{\delta}$ subset of a compact metrizable space.
This shows that a second countable, locally compact space is Polish. 

The Baire Category Theorem also applies to general locally compact spaces, 
metrizable or not, and so they are non-meager as well.

\medskip


Under the blanket assumption that $X$ is Hausdorff, and 
writing ``c.d.b., n-m.'' for ``countable density basis, non-meager," 
the system of hypotheses and implications 
for a topological dynamic system $(X,f)$ looks like this:

\medskip

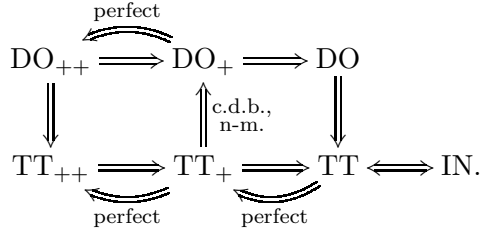
\centerline{
\xymatrix{
\mr{DO}_{++}  	\ar@{=>}[r] 	\ar@{=>}[d] 						& 
\mr{DO}_+  	\ar@{=>}[r]  	\ar@/_1pc/@{=>}[l]_{\textrm{perfect}}	&
\mr{DO}  					\ar@{=>}[d] 						&	\\
\mr{TT}_{++}	\ar@{=>}[r] 									&
\mr{TT}_+ 		\ar@{=>}[r]	\ar@/^1pc/@{=>}[l]^{\textrm{perfect}}
			\ar@{=>}[u]_{\substack{ \mr{c.d.b.,}\\\mr{n{\mhyphen}m.} } }&
\mr{TT} 		\ar@{<=>}[r] 	\ar@/^1pc/@{=>}[l]^{\textrm{ perfect}} 	&
\mr{IN}.
}
}

\medskip


Theorem \ref{theo} is thus proven:
\begin{proof}
(a) The equivalence between $\mr{TT}$ and $\textrm{IN}$ is Proposition \ref{INTT}.
If $X$ is perfect, $\mr{TT}_{++} \iff \mr{TT}_+ \iff \mr{TT}$
by the implications (\ref{impl1}) and by
Proposition \ref{prop1}.

(b) This is Proposition \ref{prop3a}.

(c) This is Proposition \ref{prop5}.
\end{proof}

\section{The case of isolated points}

\nd\emph{N.B. In this section, our standing assumptions are that 
$X$ is Hausdorff with at least one isolated point
and $(X,f)$ is topologically transitive.}

\medskip

In this section we describe what happens when there are isolated points 
and in the process obtain the following further equivalences.

\begin{theo}\label{theo2} 
Let $(X,f)$ be a dynamical system. 
Assume that $X$ contains at least one isolated point.
\begin{enumerate}
\item[(a)] If $(X,f)$ satisfies $\mr{TT}_+$ then $X$ is finite and consists of a single
periodic orbit.  
In particular, $\mr{TT}_{+}$, $\mr{TT}_{++}$, and $\mr{DO}_{++}$ are equivalent
for dynamical systems with isolated points, 
and each implies all the other conditions.

\item[(b)] If $(X,f)$ is topologically transitive then the set of isolated points
is contained in a single orbit which is dense in $X$.  
In particular, $\mr{TT}$ is equivalent to $\mr{DO}$
for dynamical systems with isolated points.\end{enumerate}
\end{theo}

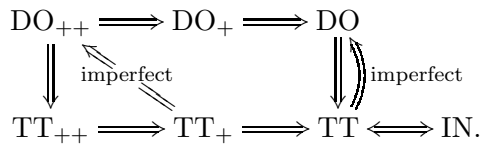
\centerline{
\xymatrix{
\mr{DO}_{++}  	\ar@{=>}[r] 	\ar@{=>}[d] 						& 
\mr{DO}_+  	\ar@{=>}[r]  									&
\mr{DO}  		\ar@{=>}[d] 									&	\\
\mr{TT}_{++}	\ar@{=>}[r] 									&
\mr{TT}_+ 		\ar@{=>}[r]	\ar@{=>}[ul]|{\textrm{imperfect}}		&
\mr{TT} 		\ar@{<=>}[r] 	\ar@/_.75pc/@{=>}[u]_{\textrm{imperfect}} 	&
\mr{IN}.
}
}

\bigskip

Combining this result with Theorem \ref{theo} thus yields further results.

\begin{cor}
Let $(X,f)$ be a dynamical system.
$\mr{TT}_+$ and $\mr{TT}_{++}$ are equivalent.
If $X$ is non-meager with countable density basis, 
$\mr{TT}_+$ implies $\mr{DO}_+$.
\end{cor}

We write $\Iso$ for the set of isolated points of $X$.
For a point $x \in X$ we will, somewhat abusively, let $x$ also
denote the singleton $\{ x \}$, allowing context to determine the reference.

It is easy to describe the case of a transitive homeomorphism.

\begin{prop}\label{IsoHomeo}
Let $f$ be a topologically transitive homeomorphism and $x \in X$ an isolated point.
The orbit $O_{\pm}(x)$ is exactly $\Iso$ and it is dense in $X$.
$f$ is point transitive iff $X$ is finite, 
in which case it consists of a single periodic orbit.
\end{prop}

\begin{proof}
All points of $O_{\pm}(x)$ are isolated since $f$ is a homeomorphism. 
Since $O_\pm(x)$ is invariant and opene, by Proposition \ref{InvOpene}(c),
it is dense. In particular, for $y \in \Iso$, it contains $y$.
Thus $O_\pm(x) = \Iso$.

If the points of the orbit sequence are all distinct, 
then clearly none of them are transitive points and the orbit is infinite.
Since $X \setminus O_{\pm}(x)$ is a proper, closed, invariant subset of $X$,
none of these points can be transitive points either.

On the other hand,
if $f^j(x) = f^k(x)$ for some $j > k$,
then $x = f^{j-k}(x)$,
so $x$ is periodic.
The orbit $O(x) = O_{\pm}(x)$ is closed because it is finite. 
As it is dense,  there are no other points in $X$,
and all points are transitive points.
\end{proof}

When $f$ is merely a continuous map, there are a few more cases.
We now describe the structure of the isolated point set in a series of lemmas. 
Recall that in general $ O_\pm(x) = \bigcup_{k \in \Z} f^k(x)$ 
is the union of all of the orbit sequences through $x$.

\begin{lem}\label{lem8a}
If $x,y \in \Iso$, then
$x \in O(y)$ or $y \in O(x)$.
Thus for any $x \in \Iso$ the set $ O_\pm(x)$ contains $\Iso$.
\end{lem}

\begin{proof}
$x$ and $y$ are open, so by $\mr{TT}$ we have
$N_+(x,y) \cup N_+(x,y) \neq \0$.
\end{proof}

\begin{lem}\label{lem8b}
If $x \in \Iso$ and $f^{-1}(x)$ contains more than one point,
then $x$ is periodic, and $f^{-1}(x)$ consists of exactly two points
(one of which lies in $O(x)$ by Lemma \ref{lem8a}).
\end{lem}

\begin{proof}
Let $U, V$ be disjoint opene subsets of $ f^{-1}(x)$
labeled so that $N_+(U,V) \neq \0$. 
Let $k$ be the smallest element of $N_+(U,V)$; 
since $U$ and $V$ are disjoint, $k > 0$.
Let $y \in U$ such that $f^k(y) \in V$.
Because $U \cup V \subset f^{-1}(x)$,
we know $x = f(y)$ and
$x = f\big(f^{k}(y)\big) = f^k\big(f(y)\big) = f^k(x)$.
Thus $x$ is a periodic point and its forward orbit $O(x) = \{ x,\ldots,f^{k-1}(x) \}$
meets $V$ at the point $f^k(y) = f^{k-1}(x)$. 
The period of $x$ cannot be a proper divisor $\ell$ of $k$, for then we would have 
$f^{\ell - 1}(x) = f^{\ell - 1}\big(f^{k - \ell}(x)\big) = f^{k-1}(x) \in V$,
so that $f^\ell(y) \in V$, contradicting minimality of $k \in N_+(U,V)$.

Suppose, for a contradiction, there exists a third point in $f^{-1}(x)$.
By shrinking $U$ and $V$ if necessary we may find a third opene
$W \subset f^{-1}(x)$ disjoint from $U$ and $V$.
By topological transitivity there exists $m \in N_+(U,W) \cup N_+(W,U)$.
If $m \in N_+(U,W)$, then $m > 0$ by disjointness
and there exists $z \in U$ such that $f^m(z) \in W$.
Again, $U \cup W \subset f^{-1}(x)$
implies $x = f\big(f^m(z)\big) = f^m\big(f(z)\big) = f^m(x)$. 
Hence $m$ is a multiple of $k$.
Thus $f^m(z) = f^{m-1}(x) = f^{k-1}(x)$ is in $W \cap V$, contradicting disjointness.
We get a similar contradiction if $m \in N_+(W,U)$.
\end{proof}

\begin{lem}\label{lem8c}
If $x$ is an isolated point,
then the preimage $f^{-1}(x)$ is a finite open set
and so consists entirely of isolated points.
\end{lem}

\begin{proof}
The preimage $f^{-1}(x)$ of an isolated point $x$ is an open set
that by Lemma \ref{lem8b} has cardinality 0, 1, or 2.
As this set is finite and $X$ is Hausdorff, the points are isolated.
\end{proof}

\begin{lem}\label{lem8d}
For any $x \in \Iso$, we have $O_-(x) \subset \Iso$, 
so that $\Iso$ is an open, $-$invariant subset of $X$.
\end{lem}

\begin{proof}
This follows from Lemma \ref{lem8c} by induction and taking unions.
\end{proof}

\begin{lem}\label{lem8e}
If $x$ is a periodic isolated point,
then $O(x) \subset \Iso$.
If, further, $O(x) \neq X$,
there is an isolated point $y$ such that
$f^{-1}\big(O(x)\big) = O(x) \cup \{y\}$.
\end{lem}

\begin{proof}
Since $x$ is periodic and isolated, 
by Lemma \ref{lem8d}, 
$O(x) \subset O_-(x) \subset \Iso$. 
If the clopen set $A = O(x)$ is a proper subset of $X$, then by Lemma \ref{lem1xxx}
the open set $f^{-1}(A) \setminus A$ is a singleton.
\end{proof}

\begin{lem}\label{lem8f}
There is at most one $x \in \Iso$
such that $f^{-1}(x)$ contains more than one point.
\end{lem}

 \begin{proof}
If $x,x' \in \Iso$,
by Lemma \ref{lem8a} we may relabel them so that $x \in O(x')$.
If their preimages each contain two points, then
$x$ and $x'$ are both periodic, by Lemma \ref{lem8b}, 
so $O(x) = O(x')$.
Then by Lemma \ref{lem8e}, 
the preimage of this orbit is $O(x) \cup \{y\}$ for some $y$.
Since $f\-(x)$ and $f\-(x')$ can each only have one point in $O(x)$,
it follows that $y \in f^{-1}(x) \cap f^{-1}(x')$,
so that $x = f(y) = x'$.
\end{proof}

%

We now can describe how $\Iso$ sits in $X$.

\begin{prop}\label{dense}
For  $x, y \in \Iso$, we have
 $O_\pm(x) = O_{\pm}(y)$
and this set is dense. 
In particular, $\Iso \subset O_{\pm}(x)$.

For $x \in \Iso$ the following conditions are equivalent:
\begin{itemize}
\item[(i)] $O(x) \subset \Iso$.
\item[(ii)] $O_{\pm}(x) = \Iso$.
\item[(iii)] $\Iso$ is $+$invariant.
\end{itemize}
These conditions imply that $\Iso$ is dense in $X$.
\end{prop}

\begin{proof}
By (\ref{NUVOpmV}), 
$\mr{TT}$ holds just if $O_\pm(U)$ is dense for each opene $U$,
so for each $x \in \Iso$, we see $O_\pm(x)$ is dense.
From Lemmas \ref{lem8d} and \ref{lem8a},
we have $y \in O_-(y) \subset \Iso \subset O_\pm(x)$.
As $O_\pm(x)$ is $+$invariant, we see $O(y) \subset O_\pm(x)$,
so $O_{\pm}(y) = O_-(y) \cup O(y) \subset O_\pm(x)$.
By symmetry they are equal.

(i) $\implies$ (ii): If $O(x) \subset \Iso$, 
then since $O_-(x) \subset \Iso \subset O_\pm(x)$
by Lemmas \ref{lem8d} and \ref{lem8a},
we have $\Iso = O_\pm(x)$.

(ii) $\implies$ (iii): $O_\pm(x)$ is $+$invariant.

(iii) $\implies$ (i): If $x \in \Iso$ and $\Iso$ is $+$invariant then $O(x) \subset \Iso$.

Since $O_\pm(x)$ is dense, $\Iso = O_\pm(x)$ implies
$\Iso$ is dense.
\end{proof}

\comment{
If $x$ is isolated and $U$ is opene, then $N(x,U) \neq \0$.
That there exists $k \in N(x,U)$ means exactly
that $O_\pm(x) = \bigcup_{k \in \Z}f^k(x)$ meets $U$.
So $O_\pm(x)$ is dense in $X$.

By Lemma \ref{lem8d},
\begin{equation}\label{denseeq}
\bigcup_{k \in \N}f^{-k}(x) \quad \subset \quad \Iso.
\end{equation}
If $y$ is also isolated then by Lemma \ref{lem8a} $\Iso \subset O_\pm(y)$ and so $x \in O_\pm(y)$.
Since $O_{\pm}(y)$ is clearly $+$invariant, $O(x) \subset O_\pm(y)$. From (\ref{denseeq})
it follows that $O_\pm(x) \subset O_\pm(y)$.
By symmetry they are equal.

(i) $\implies$ (ii): If $O(x) \subset \Iso$ then $O_\pm(x) \subset \Iso$ by (\ref{denseeq}).
So the two sets are equal by Lemma \ref{lem8a}.

(ii) $\implies$ (iii): $O_\pm(x)$ is $+$invariant.

(iii) $\implies$ (i): If $x \in \Iso$ and $\Iso$ is $+$invariant then $O(x) \subset \Iso$.

Since $O_\pm(x)$ is dense, $\Iso = O_\pm(x)$ implies
$\Iso$ is dense.
}

Putting these facts together gives a complete description 
of how isolated points occur in a topologically transitive system.


\begin{prop}\label{Iso} Assume that $(X,f)$ is topologically transitive 
and that $X$ contains isolated points. 
The system has a dense orbit.
Exactly one of the following cases occurs:
\begin{itemize}

\item[1.] There exists a unique $x \in \Iso$ such that $f^{-1}(x) = \0$. 
$\Trans = \{ x \}$, and $f(X)$ is not dense in $X$. 
In this case $\mr{DO}_+$ holds, 
but $\mr{TT}_+$ and hence $\mr{TT}_{++}$ and $\mr{DO}_{++}$ do not.
Exactly one of the following occurs:
\begin{enumerate}

\item[1a.] \emph{``$\N$'':}
$\Iso  = O_\pm(x) = O(x)$ consists of infinitely many distinct points 
in a single forward orbit, and $\Iso$ is dense in $X$.

\item[1b.] \emph{``Finite figure \textsf{9}'':}
$\Iso = O_\pm(x) = O(x) = X$ is a finite, pre-periodic forward orbit of period $\ell$.
$y = f^k(x)$ is periodic for some minimum $k > 0$, 
and $f^{-1}(y) = \{f^{k-1}(x),f^{\ell -1}(y)\}$.

\item[1c.] \emph{``$n$'':}
$\Iso = \{f^{k}(x) : 0 \leq k \leq n-1\}$ is a finite sequence of distinct points,
for some $n \geq 0$. 
For $k \geq n$, $f^k(x)$ is not isolated and $X$ is infinite.
The finite set $\Iso$ is not dense in $X$.  $O(x) = O_\pm(x)$ is dense in $X$.
\end{enumerate}

\item[2.] For every point $z \in \Iso$ the set $f^{-1}(z) \neq \0$, 
and there exists $x \in \Iso$ such that $f^{-1}(x)$ contains two points. 
In that case the point $x$ is unique and we have
\begin{enumerate}

\item[2.] \emph{``Infinite figure \textsf{9}'':} 
$\Iso = O_\pm(x)$ is an infinite, pre-periodic orbit of period $\ell$, 
and $f^{-1}(x) = \{y, f^{\ell -1}(x) \}$ for some $y$.
For all $k \in \N$, $f^{-k}(y)$ is a single isolated point.
$\Trans = \0$, and
$\Iso$ and $f(X)$ are dense in $X$.
\end{enumerate}

\item[3.] For every $x \in \Iso$, the preimage is a singleton,
and exactly one of the following occurs:

\begin{enumerate}
\item[3a.] \emph{``Figure \textsf{0}'':} $X = f(X) = \Trans = \Iso$ 
is a single periodic orbit. 
This is the only case satisfying $\mr{TT}_{++}$, $\mr{TT}_+$, or $\mr{DO}_{++}$.

\item[3b.] \emph{``$\Z$'':} 
For each $x \in \Iso$, 
the bi-infinite orbit $O_\pm(x) = \Iso$.
$\Trans = \0$,
$\Iso$ is dense in $X$,
and $f(X)$ is dense in $X$.

\item[3c.] \emph{``$-\N$'':} 
There is a unique $y \in \Iso$ such that $f(y) \not\in \Iso$.
$\Iso = \{f^{-k}(y) : k \in \N\}$ forms an infinite sequence ending in $y$ and
for $k > 0$, $f^k(y)$ is not isolated.
$\Trans = \0$
and $f(X)$ is dense in $X$. $\Iso$ may or may not be dense in $X$.
\end{enumerate}
\end{itemize}
\end{prop}


\begin{proof}
By Proposition \ref{dense}, 
$x \in \Iso$ implies $\Iso \subset O_\pm(x)$, and the latter set is dense. 

By Corollary \ref{cor1xxx}, 
$f(X)$ fails to be dense exactly when there is an isolated point $x$ --- 
necessarily unique --- such that $f^{-1}(x) = \0$, i.e. in Case 1. 

Suppose that there is some $y \in \Iso$ such that $f(y)$ is not isolated.
Then no iterate $f^k(y)$ for $k > 0$ is isolated,
since by Lemma \ref{lem8d}, $\Iso$ is $-$invariant.
Since $\Iso \subset O_\pm(y)$ we see $\Iso =O_-(y)$.
If for any $x \in \Iso$ the preimage $f^{-1}(x)$ contained more than one point,
by Lemma \ref{lem8b} we would have $x$ and hence $y$ periodic.
But then, by Lemma \ref{lem8e}, $O(y)$ would consist of isolated points,
contrary to assumption.
Thus for each point $x \in \Iso$,
the cardinality of $f^{-1}(x)$ is zero or one. 
If it is always one, we are in Case 3c, the $-\N$-shape.
If it is zero for some $x$, we are in Case 1c, the $n$-shape; 
in this case, since $\Iso$ is a proper finite subset of $X$, it is not dense.

Now suppose there is no $y \in \Iso$ such that $f(y)$ is not isolated.
Then $\Iso$ is $+$invariant, 
and so by Proposition \ref{dense}, 
$\Iso = O_\pm(x)$ for $x \in \Iso$ and so is dense. 
Now we consider these cases.

Suppose there is a periodic isolated point $x$. 
If $O(x) = X$ then we are in Case 3a, the figure \textsf{0}.
Otherwise, by Lemma \ref{lem8e}, 
$f^{-1}\big(O(x)\big) \setminus O(x)$ is a single isolated point $y$ 
and $\Iso = O_\pm(y) =
O(x) \cup O_-(y)$. 
By Lemma \ref{lem8f}, the preimage of any point of $ \Iso \setminus O(x)$,
consists of one point or none.
If  every preimage of every isolated point other than $f(y)$ is a single point,
we are in Case 2, the infinite figure \textsf{9}.
Otherwise, we are in Case 1b, the finite figure \textsf{9}.

Now suppose that $\Iso$ is $+$invariant and there are no periodic isolated points.
Then by Lemma \ref{lem8b}, 
every preimage $f^{-1}(y)$ of an isolated point $y$ contains one or zero points. 
If it is always one, we are in Case 3b, the $\Z$-shape.
Otherwise we are in in Case 1a, the $\N$-shape.

This exhausts the possibilities.

The results about transitive points are easy to check from the observation,
following from Lemma \ref{lem8a},
that
\[
x \in \Trans \implies \Iso \subset O(x).
\]
From this it follows that
(1) if $x \in \Iso$ and $f^{-1}(x) \neq \0$,
then $x \notin \Trans$ unless $x$ is a periodic point and $X = O(x)$, and
(2) if $x \notin \Iso$ then $x \notin \Trans$.
\end{proof}

We establish Theorem \ref{theo2}.
\begin{proof}
(a) If $X$ is a single periodic orbit then $(X,f)$ satisfies 
$\mr{DO}_{++}$, $\mr{TT}_{++}$, and $T_{+}$;
the remaining isolated point cases satisfy none of them.

\nd (b) 
In the isolated point cases, $\mr{TT}$ implies $\mr{DO}$
by Lemma \ref{lem8a},
and the converse as always true.
\end{proof}

We sketch what examples of these cases look like. 
Let $(\Z,s)$ be the successor dynamical system on the
discrete integers with $s(n) = n+1$.

\begin{itemize}
\item ``$\Z$'': Every example contains a copy of $(\Z,s)$ 
with $\Z$ a dense, open, invariant subset. 
Each compact example maps onto the system obtained by extending $s$ 
to the one-point compactification of $\Z$.


\item ``$\N$'': $\N$ is a $+$invariant subset of $\Z$ 
and every example contains a copy of the subsystem
$(\N,s \rest \N)$ with $\N$ as a dense, open, invariant subset. 
Each compact example maps onto the system obtained 
by extending $s$ to the one-point compactification of $\N$.

\item ``Figure \textsf{0}'':
Every example is isomorphic with $(\Z/n\Z, \bar s)$ for some 
$n > 0$, 
where $\bar s$ is induced by the successor function on 
the quotient group of integers modulo $n$.

\item ``Infinite figure \textsf{9}'':  With $n > 0$, 
define $\Z/n \N$ to be the quotient space obtained by identifying
two nonnegative integers if they are congruent modulo $n$. 
The translation $s$ on $\Z$ induces a map $\bar s$ on $\Z/n \N$ such that
the preimage of the class $\{0\}$ comprises the two classes $\{-1\}$ and $[n-1]$. 
Every example
contains a copy of this system with $\Z/n \N$ as a dense, open, invariant subset. 
Each compact example maps onto the system obtained by 
extending $\bar s$ to the one-point compactification of $\Z/n \N$.

\item ``Finite figure \textsf{9}'': With $n, k > 0$, 
the image of the set of integers at least $-k$ is a $+$invariant
subset of $\Z/n \N$, 
and every example is isomorphic to the subsystem of $(\Z/n \N,\bar s)$ 
induced by such a set.

\item ``$n$'': Let $(Y,g)$ be any point transitive system on a perfect space 
with a transitive point $y_0$. 
Let $X$ be the disjoint union of $Y$ with $\{1, \ldots, n \}$ 
and define $(X,f)$ to extend the subsystem $(Y,g)$ by
$f(i) = i+1$ for $1 \leq i < n$ and $f(n) = y_0$.

\item ``$-\N$'': This is the only interesting case. 
We describe the metric space examples.
Let $(Y,g)$ be a dynamical system with the perfect space $Y$ 
the union of two closed, $+$invariant subspaces, $Y = Y_1 \cup Y_2$. 
Assume that the subsystem $(Y_1,g)$ contains a transitive point $y_0$
and $Y_2$ contains a dense set $\{ y_{-1}, y_{-2}, \ldots \}$ satisfying
$\lim_{k \to \infty} d\big(y_{-k},g(y_{-(k+1)})\big) = 0$.
In $Y \times [0,1]$, 
let $X = \big(Y \x \{0\}\big) \cup \big\{ \big(y_{-k},\frac 1 k\big) : k = 1, 2, \ldots \big\}$.
Define $f\colon X \to X$ by
$\begin{cases}
f(y,0) = \big(g(y),0\big),								& y \in Y;		\\
f\big(y_{-(k+1)},\frac 1 {k+1}\big) = \big(y_{-k},\frac 1 k\big),	& k \geq 1;	\\ 
f\big(y_{-1},1\big) = (y_0,0). 							&
\end{cases}$

\nd The closure of $\Iso$ is $\Iso \cup (Y_2 \times \{0\})$.
Thus $\Iso$ is dense iff $Y_1 \subset Y_2 = Y$.
\end{itemize}

\section{Minimality}

\begin{df}
A dynamical system $(X,f)$ is called \emph{minimal} 
when every point is a transitive point, i.e., when $\Trans = X$.
A $+$invariant subset $A$ of $X$ is called a \emph{minimal subset}
when the subsystem $(A, f \rest A)$ is minimal.\end{df}

By definition, a minimal system satisfies $\mr{DO}_{+}$ 
and so $\mr{DO}$ and $\mr{TT}$ as well. 
As we will now see, 
it satisfies $\mr{DO}_{++}$ and so all seven conditions for topological transitivity.

Clearly,
if $X$ consists of a single periodic orbit,
then the system is minimal.
The term ``minimal" comes from parts (a) and (d) of the following proposition.

\begin{prop}\label{prop3.2}
Let $(X,f)$ be a dynamical system which is not a single periodic orbit.
\begin{enumerate}
\item[(a)] $(X,f)$ is minimal iff $X$ contains no proper, closed, $+$invariant subsets.

\item[(b)] If $(X,f)$ is minimal
then $X$ is perfect, 
$f(X)$ is dense in $X$, and
$\omega f(x) = X$ for every $x \in X$.

\item[(c)] If $X$ is compact then $X$ contains a minimal, closed, invariant subset.

\item[(d)] If $X$ is compact then $(X,f)$ is minimal iff 
$X$ contains no proper, closed, invariant subsets.

\item[(e)] If $X$ is compact and $(X,f)$ is minimal with $f$ a homeomorphism. 
then $(X,f^{-1})$ is minimal.
\end{enumerate}
\end{prop}
\begin{proof}
(a)
If $A$ is a proper, closed $+$invariant subset of $X$,
then no point $x \in A$ is transitive
since the closure of $O(x)$ is contained in $A$.
On the other hand, for any $x \in X$,
the closure $\overline{O(x)}$ is nonempty, closed, and $+$invariant, so
if $X$ contains no proper, closed, $+$invariant subsets, then $\overline{O(x)} = X$.
Thus all points are transitive. 

(b)
When the periodic orbit case is excluded,
the remaining cases of Proposition \ref{Iso} are not minimal.
In fact, each such has at most one transitive point.
Thus $X$ is perfect.
By Corollary \ref{cor1xxx},
$f(X)$ is dense in $X$, 
and by Proposition \ref{prop3a},
$X = \omega f(x)$ for all $x$.

(c)
Any nested chain of nonempty, $+$invariant closed sets $C_\alpha$ in $X$
has a closed intersection, nonempty by compactness 
and $+$invariant as well since 
\[f\big(\bigcap C_\alpha\big) \subset 
\bigcap f(C_\alpha) \subset 
\bigcap C_\alpha.\]
Thus Zorn's Lemma gives 
a minimal nonempty, closed, $+$invariant $A \subset X$. 
By compactness, $f(A)$ is a nonempty, compact, $+$invariant subset of $A$
and so equals $A$ by minimality.
That is, $A$ is invariant.
\comment{
Using Zorn's Lemma,
one can show that $X$ contains a nonempty, closed, $+$invariant subset $A$
which is minimal with respect to inclusion.
An intersection of a family of nonempty compacta which is linearly ordered
by inclusion has a nonempty intersection.
By compactness, $f(A)$ is a nonempty, closed, $+$invariant subset of $A$
and so it equals $A$.
That is, $A$ is invariant.
}

(d)
If $A$ is a proper, closed, $+$invariant subset of $X$,
then by compactness,
$\bigcap_{k \geq 0} f^k(A)$ is a proper, closed, invariant subset. 
Now use (a).

(e)
If $f$ is a homeomorphism 
then a proper, closed  $A \subset X $ is invariant for $f$ iff it is invariant for $f^{-1}$, 
so the result follows from (d).
\end{proof}

As was mentioned earlier, 
the most useful class to which Proposition \ref{prop5} applies 
is that of Polish spaces.

\begin{prop}\label{prop3.3}
Let $(X,f)$ be a point transitive system on a perfect Polish space.
With the relative topology,
the $+$invariant subset $\Trans$ is Polish
and the subsystem $(\Trans, f \rest \Trans)$ is minimal.
If $A$ is any nonempty $+$invariant subset
of $\Trans$ then $A$ is second countable and perfect
and the subsystem $(A, f \rest A)$ is minimal. \end{prop}

\begin{proof}
$\Trans$ is a dense, $G_{\delta}$ subset of $X$
by Proposition \ref{prop5} and so is Polish.
Any nonempty $+$invariant subset $A$ of $\Trans$ 
is dense because $O(x)$ is dense for any transitive point $x$.
By Proposition \ref{prop3.1}, 
such an $A$ is perfect and every point is transitive for $(A, f \rest A)$. 
As a subset of a second countable space, $A$ is second countable.
\end{proof}

\section{Examples}
In this section we show that the implications proven in Section 4 fail to go through
under weaker hypotheses. As usual, all spaces are Hausdorff.

\begin{ex}\label{ex1}
For $X$ with isolated points,
$\mr{DO}_+$ $\nimplies$ $\mr{TT}_+$, $\mr{TT}_{++}$, or $\mr{DO}_{++}$,
even if $X$ is compact metrizable.
$\Iso$ need not be $+$invariant.
\end{ex}
\begin{proof}
All cases of Proposition \ref{Iso} with an initial point, Cases 1a, 1b, and 1c,
satisfy $\mr{DO}_{+}$, with the initial point $x$ the unique transitive point.
However $U = X \setminus x$ is opene with $N_+(U,x) = \0$,
so $\mr{TT}_+$ does not hold, 
and hence not $\mr{TT}_{++}$ or $\mr{DO}_{++}$ either.
In Case 1c, $\Iso$ is not $+$invariant.
\end{proof}

\begin{ex}\label{ex2}
If $X$ is not perfect,
$\mr{DO} \nimplies \mr{DO}_{+}$ or $\mr{TT}_+$,
even if $X$ compact and metrizable and $f$ a homeomorphism.
\end{ex}
\begin{proof}
All cases of Proposition \ref{Iso} satisfy $\mr{DO}$.
Cases 2, 3b, and 3c do not satisfy $\mr{DO}_+$ or $\mr{TT}_+$.
The homeomorphism examples are all in Case 3b.
\end{proof}

We will build many of our examples by using the following:

\begin{prop}\label{prop3.1}
Let $(X,f)$ be a dynamical system and $A$ a dense, $+$invariant subset of $X$.
$(X,f)$ satisfies $\mr{TT}_{++}$, $\mr{TT}_{+}$, or $\mr{TT}$
iff the subsystem $(A,f \rest A)$ satisfies the corresponding property.
If $x \in A$, then $x$ is a transitive point for $(X,f)$
iff it is a transitive point for $(A,f \rest A)$,
and $\omega(f \rest A)(x) = A$ iff $\omega f(x) = X$.
$X$ is perfect iff $A$ is perfect.
\end{prop}
\begin{proof}
Since $A$ is dense, $U \cap f^{-k}(V) $ is opene 
iff $U \cap f^{-k}(V) \cap A$ is opene in $A$.
Since $A$ is $+$invariant,
\[
U \cap f^{-k}(V) \cap A =
(U \cap A) \cap (f \rest A)^{-k}(V \cap A)
\]
for all $k \in \N$.
Hence $N_+(U,V)$ for $f$
equals $N_+(U \cap A, V \cap A)$ for $f \rest A$.
By definition of the subspace topology,
the open sets in $A$ are exactly of the form $U \cap A$ with $U$ open in $X$.
Thus $\mr{TT}_{++}$, $\mr{TT}_+$, or $\mr{TT}$ holds for $(X,f)$ 
just if it does for $(A, f \rest A)$.
If $x \in A$ is such that $O(x)$ is dense in $A$ (or $\omega (f\rest A)(x) = A$),
then $O(x)$ is dense in $X$ (resp. $\omega f(x) = X$).

If $x$ is an isolated point of $X$ then the open set $x$ meets the dense set $A$ 
and so $x = x \cap A$ is an isolated
point of $A$.  On the other hand, if $x = U \cap A$ is an isolated point of $A$ 
with $U$ open in $X$ then
$U \setminus x$ is an open set disjoint from $A$ and so is empty.  
That is, $U = x$ is an isolated point of $X$. It follows
that $X$ is perfect iff $A$ is.
\end{proof}


If $A$ is any compact space with more than one point,
then the product space $A^\Z$ is a perfect, compact space.
It is the space of bi-infinite sequences in $A$
and on it we define the (left) \emph{shift homeomorphism}
$s \colon A^\Z \to A^\Z$ by
\[
s(x)_i \: = x_{i+1}.
\]
Because a basic open set in $A^\Z$ restricts only finitely many coordinates,
it is easy to see that $(A^\Z,s)$ is topologically transitive.
In fact, for any opene $U, V \subset A^\Z$,
$N(U,V)$ contains all but finitely many integers.
Thus the system satisfies $\mr{TT}_{++}$.

\begin{ex} \label{ex3}
If $X$ is not separable, then $\mr{TT}_{++} \nimplies \mr{DO}$,
even if $X$ is perfect and compact
and $f$ is a homeomorphism.
\end{ex}

\begin{proof}
Write $\c = 2^{\aleph_0}$ for the cardinality of $\R$. 
If $A$ is a compact space of cardinality greater than $2^\c$
(e.g. $[0,1]^\kappa$ or $\{0,1\}^\kappa$ for $\kappa > 2^\c$), 
then $A$ cannot contain a countable dense set,
so neither can $A^\Z$.
Therefore $\mr{DO}$ does not hold of $(A^\Z,s)$ even though $\mr{TT}_{++}$ does.
%
\end{proof}

\begin{ex}\label{ex3.2a}
There exists $(X,f)$ with $X$ perfect and compact but not admitting a countable 
density basis 
--- and hence not second countable and thus not metrizable --- 
which satisfies $\mr{DO}_{++}$.
\end{ex}
\begin{proof}
We will sketch the construction.  For details see \cite{A2}, p. 91.

Let $K$ be the Cantor set and $A$ a compact space.
Via the identifications $(A^\Z)^K \iso A^{\Z \x K} \iso (A^K)^\Z$,
the shift $s$ on $(A^K)^\Z$
induces a homeomorphism $s_*$ on $(A^\Z)^K$.
In more detail, given $\alpha\colon K \to A^{\Z}$
then $s_*(\alpha) = s \circ \alpha\colon K \to A^\Z$.
Since $\big((A^K)^\Z,s\big)$ satisfies $\mr{TT}_{++}$ , 
so does $\big((A^\Z)^K,s_*\big)$.

When $A = 2 \: = \{ 0, 1 \}$,
the sequence space $2^\Z$ is perfect, compact, and metrizable, but
the compact Hausdorff space $(2^\Z)^K$ is not second countable, 
hence not metrizable.
In fact we have the following:

\begin{prop}\label{proddensebasis} 
If $Y$ contains at least two points and $K$ is uncountable then the product space
$Y^K$ does not admit a countable density basis. 
\end{prop}
\begin{proof} 
If $\D$ is a density basis for a space and $\mathcal{B}$ 
is a basis for the topology then
we can choose for each $U \in \D$ 
a nonempty $V \in \mathcal{B}$ such that $V \subset U$.  
These choices determine a family $\tilde{\D} \subset \mathcal{B}$ 
which is of cardinality at most that of $\D$.  
If a set meets every element of $\tilde{\D}$ then it meets every element of $ \D$
and so is dense.  
Thus $\tilde{\D}$ is a density basis consisting of basic open sets.
Now suppose that $ \D$ is a countable collection of basic open sets of $Y^K$, 
that is, that each is of the form $\bigcap_{i \in F} \ \pi_i^{-1}(U_i)$ 
where $F$ is a finite subset of $K$ and each $U_i \subset Y$ is open.
Taking the union of all of the index sets $F$ associated with elements of $ \D$ 
we obtain a countable $I \subset K$.  
Let $U_1, U_2$ be disjoint opene subsets of $Y$ and let $j \in K \setminus I$. 
Since the opene sets $\pi_j^{-1}(U_1),\pi_j^{-1}(U_2) \subset Y^K$ are disjoint, 
neither is dense.  
But each meets every element of $ \D$. 
Thus $ \D$ is not a density basis.
\end{proof}

%
%

Let $C(K, 2^\Z)$,
the set of continuous functions from the Cantor set $K$ to $2^\Z$.
Let $d$ be a metric on $2^\Z$ inducing the product topology,
for example $d(x,y) = \max_{n \in \Z} 2^{-|n|}|x_n - y_n|$,
and topologize $C(K, 2^\Z)$ with the \emph{sup metric} $\rho$;
that is to say, $\rho(\a,\b) = \max_{k \in K}d(\a(k),\b(k))$.
This yields the topology of uniform convergence.
This metric is complete
because the uniform limit of continuous functions is continuous.
Furthermore, because $K$ is a Cantor set,
the collection of locally constant functions with image in a countable dense subset of
$2^{\Z}$ form a countable dense subset of $C(K,2^\Z)$.
Hence  $C(K,2^\Z)$  is Polish. The inclusion
map $J \colon C(K,2^\Z) \to (2^\Z)^K$ is continuous and injective
but it is not a homeomorphism onto its image.
In particular, we cannot immediately use Proposition \ref{prop3.1}
to show that $s_*$ is topologically transitive on $C(K,2^\Z)$.
However, it is not hard to show that
the restriction of $s_*$
to the set of locally constant functions in $C(K,2^\Z)$ is topologically transitive,
and so we can use Proposition \ref{prop3.1}
to see that $s_*$ is topologically transitive on $C(K,2^\Z)$.
As this is a Polish space,
Proposition \ref{prop5} implies there exists a transitive point $\alpha$ 
for $s_* \rest C(K,2^\Z)$.
That is, the forward orbit $O(\alpha)$ is dense in $C(K,2^\Z)$.
Because $J\big(C(K,2^\Z)\big)$ is dense in $(2^\Z)^K$,
it follows that $J\big(O(\alpha)\big)$ is dense in $(2^\Z)^K$.
That is, the orbit of $\alpha$ is dense in $(2^\Z)^K$.
It follows from Proposition \ref{prop3a}
that $s_*$ satisfies $\mr{DO}_{++}$ on $(2^\Z)^K$.
\end{proof}

As noted above, when $A = 2 =\{ 0, 1 \}$,
the sequence space $2^\Z$ is perfect, compact, and metrizable.
The shift homeomorphism on sequences of the two symbols
$0$ and $1$ defines a system $(2^\Z,s)$ satisfying $\mr{DO}_+$ 
(and hence $\mr{DO}_{++}$), 
because a point is transitive whenever every finite sequence of
symbols --- every \emph{word} formed from the two-symbol alphabet ---
appears on the positive side of the sequence.
On the other hand, a sequence is periodic
exactly when it is a periodic point for the homeomorphism $s$.
It is easy to see that the set of periodic points is dense in $2^\Z$.
The point $\bar 0$ defined by $\bar 0_j = 0$ for all $j \in \Z$ is a fixed point.
Define
\[
\Transz \: =
\{ x \in \textit{Trans}_s   :   \exists N \in \Z \ \forall n < N \ (x_n = 0) \}.
\]
That is, $x \in \Transz$ when every finite word 
appears on the positive side of the sequence, but on the negative side the
values are eventually $0$, so the forward orbit $O(x)$ is dense 
but as $k \to \infty$ the sequence $s^{-k}(x)$ converges to $\bar 0$.

\begin{ex} \label{ex4}
For $X$ meager, $\mr{TT}_{++} \nimplies \mr{DO}$,
even if $X$ is perfect, second countable, and metrizable
and $f$ is a homeomorphism.
\end{ex}
\begin{proof}
Let $A$ be the set of periodic points for $s$ in $2^\Z$.\footnote{\ 
This example is actually conjugate to the set of periodic orbits of the tent map.}
This is a dense, invariant subset of $2^\Z$
and so $(A,s \rest A)$ satisfies $\mr{TT}_{++}$ by Proposition \ref{prop3.1}.
Every orbit is finite and so is nowhere dense in the perfect space $A$.
\end{proof}

If $(X,f)$ is topologically transitive and $X$ is compact and perfect,
then $f(X)$ is a compact subset of $X$,
dense by Corollary \ref{cor1xxx}, 
and so equals $X$.
That is, $f$ is surjective.

\begin{ex} \label{ex5}
If $X$ is not compact, then $\mr{DO}_{++}$ implies neither that $f$ is surjective 
nor that $\Trans$ is invariant,
even if $X$ is perfect, second countable, and locally compact.
\end{ex}
\begin{proof}
Begin with $2^\Z$ and choose $x \in \Transz$.
Let $X = 2^\Z \setminus \big(\{ \bar 0 \} \cup \{ s^{-k}(x) : k \geq 1 \}\big)$ 
and $f = s \rest X$.
Clearly, $X$ is a dense, $+$invariant subset of $2^\Z$.
Since $X$ is open, it is locally compact. 
The subsystem $(X,f)$
is topologically transitive by Proposition \ref{prop3.1},
which also implies that $X$ is perfect.
Hence $(X,f)$ satisfies $\mr{DO}_{++}$ by
Proposition \ref{prop5}  and Proposition \ref{prop3a}.
On the other hand, $x \in \Trans \setminus f(X)$. 
Thus neither $X$ nor $\Trans$ is invariant.
\end{proof}

\begin{ex}\label{ex6a}
If $X$ is meager,
$(X,f)$ minimal with $f$ a homeomorphism, 
does not imply $(X,f^{-1})$ satisfies $\mr{DO}_+$,
even if $X$ is perfect and second countable.
If $X$ is meager, 
$\mr{DO} \nimplies \mr{DO}_+$,
even if $f$ is a homeomorphism and $X$ is perfect and second countable.
\end{ex}
\begin{proof}
Begin with $(2^\Z,s)$. 
Let $f = s \rest \Transz$, a homeomorphism on the
meager invariant set $X = \Transz$.
By Proposition \ref{prop3.3},
the subsystem $(X,f)$ is minimal,
and $X$ is perfect.
However, for the inverse system $(X,f^{-1})$
no orbit is a transitive point
since the orbit sequences all converge in $2^\Z$ to $\bar 0$.
Thus the inverse system does not satisfy $\mr{DO}_+$.
Of course, it does satisfy $\mr{DO}$.
\end{proof}

Recall that for a point $x \in X$ under a map $f \colon X \to X$
we define
$O_\pm(x) = \bigcup_{k \in \Z} f^k(x)$ even when $f$ is not a homeomorphism.
If $f^{-1}(x)$ is countable for every $x \in X$ then each $O_{\pm}(x)$ is countable.

\begin{ex}\label{ex6b} 
The existence of a dense $O_{\pm}(x)$ for some point $x \in X$
(as opposed to a dense orbit sequence)
need not imply $\mr{TT}$ for $(X,f)$,
even if $X$ is a perfect, compact, metrizable space with each $f^{-1}(x)$ finite.
\end{ex}
\begin{proof}
Choose $x \in \Transz$ so that its forward $s$-orbit $O(x)$ is dense in $2^\Z$.
Now let $X_1 \: = 2 \times  2^\Z$ be a disjoint union of two copies of $2^\Z$,
and define $f_1(i,x) \: = (i,s^{-1}(x))$ for $i = 0,1$.
Define $E \subset X_1$ by
\[
E = 2 \x \big(\{\bar 0\} \cup \{ s^{-k}(x) : k \geq 1\}\big).
\]
Clearly, $E$ is closed and $+$invariant under $f_1$.
Let $\pi \colon X_1 \to X$ be the quotient space projection
identifying all of the points of $E$ together to define a single point $e$.
The homeomorphism $f_1$ induces a continuous map $f$ on $X$.
The preimage $f^{-1}(e)$ of the fixed point $e$ is $\{ e, [(0,x)], [(1,x)] \}$.
Hence $O_{\pm}(e)$ is $e$ together with $2 \times O(x)$,
which is dense.
On the other hand,
$\pi(\{0\} \times 2^\Z) \setminus \{e\}$
and $\pi(\{1\} \times 2^\Z) \setminus \{e\}$
are disjoint $-$invariant opene subsets of $X$.
Hence $(X,f)$ is not topologically transitive by Proposition \ref{InvOpene}.
\end{proof}

\begin{ex} \label{ex6}
For $f$ a homeomorphism and $X$ not compact,
$(X,f)$ minimal $\nimplies$ $(X,f^{-1})$ minimal,
even if $X$ is perfect and Polish.
\end{ex}
\begin{proof}
Begin again with $(2^\Z,s)$.
By Proposition \ref{prop3.3},
the subsystem $(\textit{Trans}_s,s \rest \textit{Trans}_s)$ is minimal,
with $\textit{Trans}_s$ perfect and Polish.
Since $\textit{Trans}_s$ is both $+$invariant and $-$invariant,
it is invariant for the homeomorphism $s$,
and so $s \rest Trans_s$ is a homeomorphism.
On the other hand,
the subset $\Transz$ consists of points
which are not transitive points for $s^{-1}$,
and so $\big(\textit{Trans}_s,(s \rest \textit{Trans}_s)^{-1}\big)$ is not minimal.
\end{proof}

We conclude with an issue which remains open as far as we know:

\begin{ques}  How far can we push Proposition \ref{prop5}?
\begin{itemize}
\item  Does there exist a topologically transitive $(X,f)$ 
with $X$ compact and separable which is not point transitive, 
or even fails to have a dense orbit sequence, i.e.. does not satisfy $\mr{DO}$?

\item Does there exist a point transitive $(X,f)$ 
with $f$ a homeomorphism and $X$ compact (and necessarily separable)
but for which $(X,f^{-1})$ is not point transitive?

%

\end{itemize}
\end{ques}
\vspace{1cm}


\section{Appendix: Density Bases}

Recall our definition:

\begin{df}\label{denbasdef} 
A \emph{density basis} $\D$ for a space $X$ 
is a collection of opene subsets of $X$
such that if $A \subset X$ meets every $U \in \D$, 
then $A$ is dense in $X$. 
\end{df}

In other words, $\D$ is a density basis just when 
in order that a set $A \subset X$ meet each opene $V \subset X$, 
it suffices $A$ meet each $U \in \D$.
Evidently each basis is a density basis, then, though the converse need not hold
(see Example \ref{LexOrder}).
In particular, a second countable space admits a countable density basis.

\begin{prop}\label{nhdBasis}
If $\D$ contains a neighborhood basis for 
each point of a dense set $D \subset X$,
then $\D$ is a density basis.
In particular, if $X$ is separable and first countable 
then $X$ admits a countable density basis.
\end{prop}
\begin{proof}
Any closed set which meets every element of $\D$ 
must then contain $D$ and so must equal $X$.  
Thus $\D$ is a density basis. 
\end{proof}

\begin{ex}\label{LexOrder}
There exists a compact Hausdorff space $X$ 
which is separable and first countable
(and hence admits a countable density basis by Proposition \ref{nhdBasis}),
but which is not second countable and so is not metrizable.
\end{ex}
\begin{proof}
If the product set $X = [0,1] \times \{ 0,1 \}$ is ordered lexicographically and
given the associated order topology, 
then $X$ becomes a compact space with uncountably many clopen sets 
and so is not second countable. 
It is separable and first countable.
\end{proof}
In addition, 
if $Y$ is a separable metric space which is not compact then 
the Stone--\v{C}ech compactification $X = \beta Y$ is not metrizable 
but does admit a countable density basis 
(see Proposition \ref{dencor4} below).
On the other hand,
if $\D$ is a countable density basis, 
then by choosing a point from each set in $\D$, 
we obtain a countable dense subset, 
and so the space is separable.

\medskip

\nd\emph{N.B. 
Throughout the rest of this appendix we assume that our spaces 
are regular as well as Hausdorff.}

\begin{prop}\label{denprop1}
Let $\D$ be a collection of opene subsets of $X$. The following are equivalent:
\begin{enumerate}
\item[(i)] $\D$ is a density basis for $X$.
\item[(ii)] If a closed $C \subset X$ meets every $U \in \D$, 
then $C = X$.
\item[(iii)] If an open $V \subset X$ meets every $U \in \D$,
then $V$ is dense in $X$.
\item[(iv)] For every opene $V \subset X$, 
the set $\bigcup \{ U \in \D : U \subset V \}$ is dense in $V$.
\item[(v)]  Every opene $V \subset X$ contains some $U \in \D$.
\end{enumerate}
\end{prop}

\begin{proof}
(v) $\implies$ (i):  
If $A$ meets every element of $\D$ then by (v) it meets every
opene set and so is dense.

(i) $\implies$ (ii): 
If $C$ is closed and dense in $X$, then $C = X$.

(ii) $\implies$ (iii): 
If $V$ meets every element of $\D$, then so does its closure $\ol{V}$. 
By (ii), then, $\ol{V} = X$.

(iii) $\implies$ (iv): 
If $V$ is opene and $W $ is an arbitrary opene subset of $V$,
then by regularity there exists a closed set $B \subset W$ with nonempty interior.  
Since $X \setminus B$ is open but not dense,
by (iii), $X \setminus B$ is disjoint from some $U \in \D$.  
That means $U \subset B \subset W$. 
It follows that the union of those $U \in \D$ contained in $V$ is dense in $V$.

(iv) $\implies$ (v): Obvious.
\end{proof}

\begin{prop} \label{denprop2} Let $D$ be a dense subset of $X$.
\begin{enumerate}
\item[(a)]  If $\D$ is a density basis for $X$ then
\[
\D \wedge D \: = \{ U \cap D : U \in  \D \}
\]
is a density basis for $D$ with the relative topology induced from $X$.

\item[(b)] If $\D$ is a density basis for $D$ with the
subspace topology inherited from $X$, then
\[
\bar{\D} \: = \{ (\overline{U})^{\circ}  : U \in  \D \}
\]
is a density basis for $X$, where the closure and interior are taken in $X$.
\end{enumerate}
\end{prop}

\begin{proof} 
(a) 
If $A \subset D$ meets every element of $\D \wedge D$, 
then it meets every element of $\D$, 
and so is dense in $X$ and hence in $D$.

(b) 
If $U \in \D$, 
then $U$ is opene in $D$, 
so there exists $G$ open in $X$ such that $U = G \cap D$. 
Since $D$ is dense and $G$ is open, 
$U$ is dense in $G$ and so $G \subset (\overline{U})^{\circ}$.  
Thus the elements of $\bar{\D}$ are all nonempty.

If $V \subset X$ is open, let $B $ be a closed subset of $X$ 
with $B \subset V$ and the interior of $B$ nonempty; 
such a $B$ exists by the regularity of $X$.
Because $D$ is dense, $B^{\circ} \cap D$ is opene in $D$.  
By Proposition \ref{denprop1}(v) there exists $U \in \D$
such that $U \subset B^{\circ} \cap D$.  
Thus $\ol U \subset \ol{B^{\o}} \cap \ol D = \ol{B^\o}$
and $(\ol{U})^{\o} \subset \big(\ol{B^\o_{\vphantom{X}}}\big)^\o = B^\o \subset V$.
By Proposition \ref{denprop1}(v) again, 
$\bar{\D}$ is a density basis for $X$.
\end{proof}

%

\nd{\bfseries Remark.} 
Let $x \in D$. By using essentially the same pair of arguments we can see that
if $\D$ is a neighborhood base in $X$ for $x$ 
then $\D \wedge D$ is a neighborhood base in $D$ for $x$, 
and if $\D$ is a neighborhood base for $x$ in $D$ then $\bar{\D}$ is a
neighborhood base for $x$ in $X$. 
Thus $x$ has a countable neighborhood base in $X$ iff it does in $D$.

\medskip

From Proposition \ref{denprop2} we immediately see

\begin{cor}\label{dencor1} 
Let $D$ be a dense subset of $X$ and give $D$ the subspace topology 
inherited from $X$.  
$X$ admits a countable density basis iff $D$ does. 
\end{cor}

\begin{cor}\label{dencor2} 
If $X$ admits a countable density basis 
then every dense subset of $X$ is
separable with respect to the subspace topology.
\end{cor}

\begin{proof}
By Corollary \ref{dencor1}, 
any dense subset $D \subset X$ also admits a countable density basis $\D$.
Choosing a point from every member of a density basis for $D$, 
we obtain a dense subset of $D$ 
by Proposition \ref{denprop1}(v).
\end{proof}

A continuous map $h\colon X \to Y$ is called \emph{irreducible} 
when each $A \subset X$ is dense in $X$ iff 
$h(A)$ is dense in $Y$. 
In particular, such an $h$ has a dense image. 
In general, if $A$ is dense in $X$, then $h(A)$ is dense in $h(X)$, 
and so is dense in $Y$ when $h$ has a dense image. 
It is the converse implication which is restrictive. 
The map $h$ is called \emph{weakly almost open} if 
$\big( \overline{ h(U) } \big)^\o \neq \0$ for every opene $U \subset X$. 
Assume the image of $h$ is dense in $Y$.
If $h$ is not weakly almost open, there is an opene $U \subset X$ 
such that $\big(\overline{h(U)}\big)^\o = \0$. 
Then $G = Y \setminus \overline{h(U)}$ is open and dense in $Y$ 
while $V = h^{-1}(G)$ is not dense because it it disjoint from $U$. 
Because $G$ is open and $h$ has dense image, 
$h(V) = h(X) \cap G$ is dense in $G$ and so in $Y$, 
and thus $h$ is not irreducible.
This shows that an irreducible map is weakly almost open.

\begin{prop}\label{denprop3}
Let $h\colon X \to Y$ be a continuous map.
\begin{enumerate}
\item[(a)] 
Assume $h$ is irreducible. 
If $\D$ is a density basis for $Y$ then
\[
h^*\D  \: =  \{ h^{-1}(U) : U \in \D \}
\]
is a density basis for $X$.  
In particular, if $Y$ admits a countable density basis, then so does $X$.

\item[(b)]
Assume $h$ is weakly almost open with dense image. 
If $\D$ is a density basis for $X$ then
\[
h_*\D \: = \left\{ \big(\overline{h(U)}\big)^{\circ} : U \in \D \right\}
\]
is a density basis for $Y$. 
In particular, if $X$ admits a countable density basis, then so does $Y$.
\end{enumerate}
\end{prop}

\begin{proof} (a) 
Assume $A \subset X$ meets each $h^{-1}(U)$ for $U \in \D$.
Then $h(A)$ meets $U$ in the density basis $\D$ 
and hence is dense in $Y$. 
Since $h$ is irreducible, $A$ is dense. 

(b) 
Because $h$ is weakly almost open, 
$\big(\overline{h(U)}\big)^{\circ} \neq \0$ for every opene $U$. 

Suppose an open $V \subset Y$ meets each element of 
$\big(\overline{h(U)}\big)^{\circ}$ of $h_*\D$.
Then it meets each $h(U)$ and so $h^{-1}(V)$ meets each $U \in \D$.
As $\D$ is a density basis, $h^{-1}(v)$ is dense in $X$.
It follows that $h\big(h^{-1}(V)\big) = V \cap h(X)$ is dense in $h(X)$. 
Since $h(X)$ is dense in $Y$ and $V \cap h(X) \subset V$,
it follows $V$ is dense in $Y$. 
That $h_*\D$ is a density basis follows from Proposition \ref{denprop1}.
\end{proof}

A compactification of $X$ is a map $j \colon X \to X^*$ 
with $X^*$ compact, 
$j(X)$ dense in $X^*$,
and $j$ a homeomorphism onto its image, endowed with the subspace topology. 
A Hausdorff space admits a compactification iff it is completely regular. 
The maximum compactification of a completely regular Hausdorff space $X$ 
is the Stone--\v{C}ech compactification $j_{\beta}:  X \to \beta X$,
for which there exists a --- necessarily unique --- continuous map 
$h\colon \beta X \to X^*$ such that $j = h \circ j_{\beta}$. 
Levy and McDowell observed that the map from $\beta X $ to $X^*$ is irreducible
(\cite{LMcD}, Lemma 2.1). 
This follows from the observation that  $h^{-1}\big(j(x)\big) = \{ j_{\beta}(x) \}$
for all $x \in X$ because $j$ is a homeomorphism onto its image. 
From Corollary \ref{dencor1} we immediately have

\begin{prop}\label{dencor4} Let $X$ be a completely regular Hausdorff space. 
The following are equivalent:
\begin{enumerate}
\item[(a)]  $X$ admits a countable density basis.
\item[(b)]  $X$ has a compactification $j\colon X \to X^*$ such that $X^*$ admits a countable density basis.
\item[(c)]  For every compactification $j\colon X \to X^*$, $X^*$ admits a countable density basis.
\item[(d)]  The Stone--\v{C}ech compactification $\beta X$ of $X$ admits a countable density basis.
\end{enumerate}
\end{prop}

Compare \cite{LMcD}, Theorem 3.1.

\section{Appendix: Non-Hausdorff Spaces}

For the sake of generality we remove the Hausdorff restriction here.

\begin{prop}\label{other} 
Let $(X,f)$ be a dynamical system.
\begin{enumerate}
\item[(a)] If $X$ is a perfect $T_1$ space,
then $\mr{DO} \implies \mr{TT}_+$.
\item[(b)]
If $f$ is surjective and $x \in \Trans$,
then $\omega f (x) = X$; i.e., $\mr{DO}_+ \implies \mr{DO}_{++}$.
\end{enumerate}
\end{prop}

	

\begin{proof}
\nd (a)
    Let opene $U,V \subset  X$
    and a dense orbit sequence $O = \ang{x_j}$ of $f$ be given;
    we want to find a $p \in \N$ such that $f^p(U) \cap V \neq \0$.
    
    By $\mr{DO}$, the sets $N \: = \{n : x_n \in U\}$ and $M \: = \{m : x_m \in V\}$
    are nonempty. If there is no $p \in \N$ such that $f^p(U) \cap V \neq \0$,
    it follows that each element of $N$ is greater than any element of $M$,
    so the maximum element $m \in M$ and the minimum $n \in N$ are defined.
    Since $m-n < 0$ and $f$ is continuous,
    $f^{m-n}(U) \cap V \ni x_m$ is open.

    Now $O \less F$ is dense for any finite $F \subset X$.%
\footnote{\ 
    Assume $D$ is dense. 
    By induction it suffices to show that 
    $D \less \{x\}$ is still dense for any $x \in D$. 
    If for some opene $U$ we had $U \cap D = \{x\}$,
    that would mean $(U\less\{x\}) \cap D = \0$, 
    where $U \less \{x\}$ is opene (since $X$ is $T_1$ and perfect),
    contradicting density.
    Thus each opene $U$ meets $D \less \{x\}$.
    }    
    As $O \setminus \{x_m\}$ is dense, it contains some element $x_k \neq x_m$ 
    of the open set $f^{m-n}(U) \cap V$;
    since $m \in M$ was maximal, evidently $k < m$. 
    Let $m_1$ be the greatest such $k$.
    Again by density of $D \setminus \{x_m,x_{m_1}\}$,
    there is a greatest $m_2 < m_1$ such that
    $x_{m_2} \in f^{m-n}(U) \cap V \setminus \{x_m,x_{m_1}\}$.
    By induction, it is clear that there are arbitrarily
    small $k$ with $x_k \in f^{m-n}(U) \cap V$;
    in particular, there is $k$ such that $2m - n - k > 0$.
    Now $x_m \in V$ and also
    $x_m = f^{m-k}(x_k) \in f^{m-k}\big(f^{m-n}(U)\big) \subset f^{2m - n - k}(U)$.
    
(b)
Let $O = O(x)$ be a dense forward orbit, so that the 
closure $\overline {O}$ is $X$.
For $n \in \mathbb{N}$,
by surjectivity and continuity,
we have
$X = f^n(X) = f^n (\overline O) \subset  \overline {f^n (O)}$,
so $\overline {f^n (O)} = \ol{O\big(f^n(x)\big)} = X$.
Since this holds for all $n$, taking the intersection we get $\omega f(x) = X$.
\end{proof}

\begin{ex}\label{ex0}
If the topological space $X$ is not Hausdorff,
then $\mr{TT}_{++}$ $\nimplies$ $\mr{DO}_+$.
\end{ex}
\begin{proof}
Let $X = \N$ with the topology $\{\0,X\} \cup \{(-\infty, x) : x \in X\}$,
and let $f(x) = x+1$ be the successor map.
(Taking the subspace $(-\infty,0]$ and reassigning $f(0) = 0$ also works.)
Since $f^{-1}(X) = X$ and $f^{-1}\big((-\infty,x)\big) = (-\infty,x-1)$, $f$ is continuous.
Since all opene sets meet, we have $\mr{TT}_{++}$ trivially and $X$ is not $T_1$.
We have $\mr{DO}$, but not $\mr{DO}_+$:
no forward orbit is dense,
since $O(x)$ fails to meet $(-\infty,y)$ for $y \leq x$.
\end{proof}

\section{Acknowledgments}
The second author would like to thank his professor Boris Hasselblatt for assigning
this problem and not letting him get out of completing it,
for carefully reading over preliminary versions of this paper and making suggestions,
and for contributing two of the examples.

\end{document}